\documentclass[aap,MSNbibl,nameyear,seceqn,dvips]{arximspdf}
\usepackage{mathrsfs}

\doi{10.1214/12-AAP868} 
\volume{23}
\issue{3}
\pubyear{2013}
\firstpage{1188}
\lastpage{1218}

\makeatletter
\newcommand{\rrvert}{\vert}
\newcommand{\llvert}{\vert}
\newtheorem{theorem}{Theorem}[section]

\newtheorem{proposition}[theorem]{Proposition}
\newtheorem{lemma}[theorem]{Lemma}

\newtheorem{theorema}{Theorem}

\newproclaim{definition}{Definition}[section]
\newproclaim{example}[definition]{Example}
\newproclaim{remark}[definition]{Remark}
\newproclaim{case}[definition]{Case}
\newproclaim{condition}[definition]{Condition}

\def\aseq#1{(13.#1)}

\newcommand{\eps}{\varepsilon}
\newcommand{\eq}{\eqref}
\newcommand{\eqref}[1]{(\ref{#1})}

\newcommand{\dw}{d_{\mathrm{W}}}
\newcommand{\dk}{d_{\mathrm{K}}}

\newcommand{\toinf}{\to\infty}

\newcommand{\Exp}{\operatorname{Exp}}
\newcommand{\IE}{\mathbb{E}}
\newcommand{\IP}{\mathbb{P}}
\newcommand{\law}{\mathscr{L}}
\newcommand{\eqlaw}{\stackrel{\mathscr{D}}{=}}
\newcommand{\IN}{\mathbb{N}}
\newcommand{\IR}{\mathbb{R}}
\newcommand{\I}{{\mathrm{I}}}

\def\mid{\vert}
\makeatother

\begin{document}
\begin{frontmatter}

\title{Degree asymptotics with rates for preferential attachment
random graphs\thanksref{T1}}
\thankstext{T1}{Supported in part by NUS Research Grant R-155-000-098-133.}
\runtitle{Preferential attachment random graphs}

\begin{aug}
\author[A]{\fnms{Erol A.} \snm{Pek\"oz}},
\author[B]{\fnms{Adrian} \snm{R\"ollin}\corref{}\ead[label=e1]{adrian.roellin@nus.edu.sg}}
\and
\author[C]{\fnms{Nathan} \snm{Ross}}
\runauthor{E. A. Pek\"oz, A. R\"ollin and N. Ross}
\affiliation{Boston University, National University of Singapore
and University~of~California, Berkeley}
\address[A]{E. A. Pek\"oz\\
School of Management\\
Boston University\\
595 Commonwealth Avenue\\
Boston, Massachusetts 02215\\
USA} 
\address[B]{A. R\"ollin\\
Department of Statistics\\
\quad and Applied Probability\\
National University of Singapore\\
6 Science Drive 2\\
Singapore 117546}
\address[C]{N. Ross\\
Department of Statistics\\
University of California\\
367 Evans Hall \#3860\\
Berkeley, California 94720-3860\\
USA}
\end{aug}

\received{\smonth{8} \syear{2011}}
\revised{\smonth{4} \syear{2012}}

%
\begin{abstract}
We provide optimal rates of convergence to the asymptotic distribution
of the
(properly
scaled) degree of a fixed vertex in two preferential attachment random graph
models. Our approach is to show that these distributions are unique
fixed points
of certain distributional transformations which allows us to obtain
rates of
convergence using a new variation of Stein's method. Despite the large
literature on these models, there is surprisingly little known about the
limiting distributions so we also provide some properties and new
representations, including an explicit expression for the densities in
terms of
the confluent hypergeometric function of the second kind.
\end{abstract}

%
\begin{keyword}[class=AMS]
\kwd{60C05}
\kwd{60F05}
\kwd{05C08}
\end{keyword}
\begin{keyword}
\kwd{Random graphs}
\kwd{preferential attachment}
\kwd{Stein's method}
\kwd{urn models}
\end{keyword}

\end{frontmatter}

\section{Introduction}\label{intr}

Preferential attachment random graphs are random\break graphs that evolve by
sequentially adding vertices and edges in a random way so that
connections to
vertices with high degree are favored. Particular versions of these
models were
proposed by \citet{Barabasi1999} as a mechanism to explain the
appearance of the so-called power law behavior observed in some real world
networks; for example, the graph derived from the world wide web by considering
webpages as vertices and hyperlinks between them as edges.

Following the publication of \citet{Barabasi1999}, there has been an
explosion of
research surrounding these (and other) random growth models. This work is
largely motivated by the idea that many real world data structures can be
captured in the language of networks [see \citet{Newman2003} for a
wide survey
from this point of view]. However, much of this work is experimental or
empirical and,
by comparison, the rigorous mathematical literature on these models is less
developed [see \citet{Durrett2007} for a recent review].

For preferential attachment models, the seminal reference in the mathematics
literature is
\citet{Bollobas2001}, in which one of the main results is a
rigorous proof that the degree of a randomly chosen vertex in a particular
family of preferential attachment random graph models converges to the
Yule--Simon distribution.
Corresponding approximation results in total variation
for this and related preferential attachment models can be found in
\citet{Pekoz2010} and \citet{Ford2009}.

Here we study the distribution of the degree of a fixed vertex in two
preferential attachment models. In model~1 we start with a graph $G_2$
with two
vertices labeled one and two with an edge directed from vertex two to vertex
one. Given graph $G_n$, graph $G_{n+1}$ is obtained by adding a vertex labeled
$n+1$ and adding a single directed edge from this new vertex to a
vertex labeled
from the set $\{1, \ldots, n\}$, where the chance that $n+1$ connects
to vertex
$i$ is proportional to the degree of vertex $i$ in $G_n$ (here and
below degree means
in-degree plus out-degree). Model~2 is one
studied in \citet{Bollobas2001} and allows for self-connecting edges.
There, we
start with a graph $G_1$ with a single vertex labeled one and with an edge
directed from vertex one to itself. Given graph $G_n$, graph $G_{n+1}$ is
obtained by adding a vertex labeled $n+1$ and adding a single directed
edge from
this new vertex to a vertex labeled from the set $\{1, \ldots, n+1\}$,
where the
chance that $n+1$ connects to vertex $i\in\{1,\ldots,n\}$ is
proportional to the
degree of vertex $i$ in $G_n$ (a loop at a vertex contributes
two to its degree) and the chance that vertex $n+1$ connects to
itself is $1/(2n+1)$.

Let $W_{n,i}$ be the degree of vertex $i$ in $G_n$ under either of the models
above. Our main result is a rate of convergence in the Kolmogorov
metric
(defined below) of $W_{n,i}/(\IE W_{n,i}^2)^{1/2}$ to its distributional
limit as $n\to\infty$. Although the literature on these models is
large, there
is
surprisingly little known about these distributions. The fact that
these limits
exist for the first model has been shown by \citet{Mori2005} and
\citet{Backhausz2011} and the same result for both models can be read from
\citet{Janson2006} by relation to a generalized P\'olya urn, although the
existing descriptions of the limits are not very explicit. A further related
result
from \citet{Pekoz2010} (and the only main result there having bearing on
our work here) is that for large $i$, the distribution of $W_{n,i}$
is approximately geometric with parameter $\sqrt{i/n}$, with the error
in the
approximation
going to zero as $i\to\infty$. Thus it can be seen that if $i/n\to0$ and
$i\to\infty$, then
the distribution of $W_{n,i}/\IE W_{n,i}$ converges to a rate one exponential;
cf. Proposition~\ref{prop3}(iii) below.

The primary tool we
use here to characterize the limits and obtain rates of convergence is
a new
distributional transformation for which the limit distributions are the unique
fixed points. This transformation allows us to develop a new variation of
Stein's method; we refer to \citet{Chen2011}, \citet{Ross2011} and
\citet{Ross2007} for
introductions to Stein's method.

To formulate our main result we first define the family of densities
%
%
\begin{equation}\qquad
\kappa_s(x) = \Gamma(s)\sqrt{\frac{2}{s\pi}}\exp\biggl({
\frac
{-x^2}{2s}} \biggr) U \biggl(s-1,{ \frac{1}{2}},{ \frac{x^2}{2s}}
\biggr)\qquad\mbox{for $x>0$, $s\geq1/2$,} \label{1}
\end{equation}
where $\Gamma(s)$ denotes the gamma function and $U(a,b,z)$ denotes the
confluent hypergeometric function of the second kind (also known as the Kummer
$U$ function) [see \citet{Abramowitz1992}, Chapter 13].
Propositions~\ref{prop1}
and~\ref{prop3} below
imply that $\kappa_s$ is indeed a density for $s\geq1/2$ and we
denote by $K_s$
the
distribution function defined by the density~$\kappa_s$. Define the Kolmogorov
distance between two cumulative distribution functions $P$ and
$Q$~as\looseness=-1
\[
\dk(P,Q) = \sup_{x} \bigl\vert P(x)-Q(x) \bigr\vert.
\]\looseness=0

%
\begin{theorem}\label{thm1}
Let $W_{n,i}$ be the degree of vertex $i$ in a preferential attachment
graph on
$n$ vertices defined above and let $b_{n,i}^2 = \IE W_{n,i}^2$. For
model~$1$ with
$2\leq i\leq n$ and some constants $c,C>0$ independent of $n$,
\[
\frac{c}{\sqrt{n}}\leq\dk\bigl(\law(W_{n,i}/b_{n,i}),K_{i-1}
\bigr)\leq\frac{C}{\sqrt{n}}.
\]
For model~$2$ with $1\leq i\leq n$ and some constants $c,C>0$
independent of $n$,
\[
\frac{c}{\sqrt{n}}\leq\dk\bigl(\law(W_{n,i}/b_{n,i}),K_{i-1/2}
\bigr)\leq\frac{C}{\sqrt{n}}.
\]
\end{theorem}

%
\begin{remark}
Using Proposition~\ref{prop3} below we see an interesting difference
in the
behavior of the two models. In model~1 the limit distribution for the
degree of
the first vertex (which by symmetry is the same as that for the second vertex)
is~$K_1$, the absolute value of a standard normal random variable, whereas
in model~2
the limit distribution for the first vertex is $K_{1/2}$, the square
root of
an
exponential random variable.
\end{remark}

%
\begin{remark}
To ease exposition we present our upper bounds as rates, but the
constants are
recoverable (although probably not practical especially for large~$i$).
\end{remark}

Theorem~\ref{thm1} will follow from a more general result derived by developing
Stein's method for the distribution $K_s$. The key ingredient to our framework
follows from observing that $K_s$ is a fixed point of a certain distributional
transformation which we will refer to as the ``$s$-transformed double size
bias'' ($s$-TDSB)
transformation, which we now describe.

Recall for a nonnegative random variable $W$ having finite mean we say $W'$
has the \textit{size
bias distribution} of $W$ if
\[
\IE\bigl\{Wf(W)\bigr\} = \IE W \IE f\bigl(W'\bigr)\vadjust{\goodbreak}
\]
for all $f$ such that $\IE\vert Wf(W)\vert<\infty$ [see \citet
{Brown2006} and
\citet{Arratia2010} for surveys and applications of size biasing]. If
in addition $W$ has finite second moment, then we will write
$W''$ to denote a random variable having the size bias distribution of $W'$.
Alternatively we say $W''$ has the \textit{double size bias
distribution} of $W$
and it is straightforward to check that
%
%
\begin{equation}
\label{2} \IE\bigl\{W^2f(W)\bigr\} = \IE W^2 \IE f
\bigl(W''\bigr).
\end{equation}
Although not used below,
it is also appropriate to say that $W''$ has the \textit{square bias
distribution}
of $W$ since \eq{2} implies that we are biasing $W$ against its square.
This terminology is used in \citet{Goldstein2007} and \citet
{Chen2011} albeit
under a different notation. Now, we have the following key definition.
%
%
\begin{definition}\label{def1} For fixed $s\geq1/2$ let $U_1$
and~$U_2$ be two
independent random variables uniformly distributed on the interval
$[0,1]$, and
let $Y$ be a Bernoulli random variable with parameter $(2s)^{-1}$ independent
of $U_1$ and~$U_2$. Define the random variable
\[
V := Y \max(U_{1},U_{2})+(1-Y)\min(U_{1},U_{2}).
\]
We say that $W^*$ has the \textit{$s$-transformed double size biased} ($s$-TDSB)
distribution of $W$, if
\[
\law\bigl(W^*\bigr) = \law\bigl( V W''\bigr),
\]
where $W''$, the double size bias of $W$, is assumed to be independent of~$V$.
\end{definition}

Our next result implies that the closer a distribution is to its $s$-TDSB
transform, the closer it is to the $K_s$ distribution.
Besides the Kolmogorov
metric we also consider the Wasserstein metric between two probability
distribution functions $P$ and $Q$, defined as
\[
\dw(P,Q)=\sup_{h:\Vert h'\Vert=1}\biggl\llvert\int h(x)\,dP(x) -\int
h(x)\,dQ(x)
\biggr\rrvert.
\]

%
\begin{theorem}\label{thm2}
Let $W$ be a nonnegative random variable with $\IE W^2 =1$ and let
$s\geq1$ or
$s=1/2$. Let $W^*$ have the $s$-TDSB distribution of $W$ and be defined on
the same probability space as $W$. Then if $s\geq1$,
%
%
\begin{equation}
\label{3} \dw\bigl(\law(W),K_s\bigr) \leq8s \biggl(s+
\frac{1}{4}+\sqrt{\frac{\pi}{2}} \biggr)\IE\bigl\vert W-W^*\bigr\vert,
\end{equation}
and, for any $\beta\geq0$,
%
%
\begin{equation}
\label{4} \dk\bigl(\law(W), K_s\bigr)\leq53s \beta+
34s^{3/2} \IP\bigl[\bigl\vert W-W^*\bigr\vert>\beta\bigr].
\end{equation}
If $s=1/2$ then
\[
\dw\bigl(\law(W),K_{1/2}\bigr)\leq2 \IE\bigl\vert W-W^*\bigr\vert,
\]
and, for any $\beta\geq0$,
\[
\dk\bigl(\law(W), K_{1/2}\bigr)\leq26\beta+8 \IP\bigl[\bigl\vert
W-W^*\bigr\vert
>\beta\bigr].
\]
\end{theorem}
%
%
\begin{remark}
As can easily be read from the work
of Section~\ref{sec1} (in particular, Propositions~\ref{prop1}
and~\ref{prop3}),
the distributions $K_s$ can roughly be
partitioned into three regions where similar behavior within the range
can be expected: $s=1/2$, $1/2<s<1$ and $s\geq1$.
The theorem only covers the first and last cases as
this is what is needed to prove Theorem~\ref{thm1}.
Analogs of the results of the theorem hold in the region $1/2<s<1$,
but we have omitted them for simplicity and brevity.
\end{remark}

%
\begin{remark}
From Lemma~\ref{lem12} below and the fact that for
$h$ with bounded derivative
\[
\IE\bigl\vert h(X)-h(Y)\bigr\vert\leq\bigl\Vert h'\bigr\Vert\dw\bigl
(\law(X),
\law(Y)\bigr),
\]
we see that for $s\geq1$ or $s=1/2$ and all
$\eps>0$,
\[
\dk\bigl(\law(W), K_s\bigr)\leq\frac{\dw(\law(W), K_s)}{\eps}+\sqrt
{2}\eps.
\]
Choosing
$\eps=2^{-3/4} \sqrt{\dw(\law(W), K_s)}$
yields
\[
\dk\bigl(\law(W), K_s\bigr)\leq2^{1/4} \sqrt{\dw\bigl(
\law(W), K_s\bigr)}.
\]
Thus we can obtain bounds in the Kolmogorov metric if $\vert
W-W^*\vert$ is
appropriately bounded with
high probability or in expectation.
\end{remark}
%
%
\begin{remark}\label{rem1}
It follows from Lemmas~\ref{lem5} and~\ref{lem15} below that
$\law(W)=\law(W^*)$ if and only if $W\sim K_s$.
In the case that $s=1$, $V$ is uniform on $(0,1)$ and Proposition~\ref{prop3}
below implies that $K_1$
is distributed as the
absolute value of a standard normal random variable.
Thus we obtain
the interesting fact
that $\law(W)=\law(U W'')$ for $U$ uniform $(0,1)$ and
independent of $W''$ if and only if $W$ is distributed as the
absolute value of a standard normal variable.
This fact
can also be read from its analog for
the standard normal distribution [\citet{Chen2011}, Proposition 2.3]:
$\law(W)=\law(U_0 \vert W\vert'')$ for $U_0$ uniform $(-1,1)$ and
independent of $\vert W\vert''$ if and only if $W$
has the standard normal distribution
[see also \citet{Pitman2012}].
\end{remark}

Although there are general formulations for developing Stein's method machinery
for a given distribution [see \citet{Reinert2005}], our framework below
does not
adhere to any of these directly since the characterizing operator we
use is a
second order differential operator [see \eq{23} and \eq{25} below].
For the
distribution $K_s$, the usual first order Stein operator derived from the
density approach of \citet{Reinert2005} [following \citet{Stein1986}]
is a
complicated expression involving special functions. However, by
composing this
more canonical operator with an appropriate first order operator, we
are able to
derive a second order Stein operator (see Lemma \ref{lem7} below)
which has a form that
is amenable to our analysis. This strategy may be useful for other
distributions which have first order operators that are difficult to handle.

The usual approach to developing Stein's method is to decide on the
distribution of interest, find a corresponding Stein operator and
then derive couplings from it. The operator we use here was suggested
by the
$s$-TDSB transform which in turn arose from the discovery of a close
coupling in
the preferential attachment application. We believe this approach of using
couplings to suggest a Stein
operator is a potentially fruitful new strategy for extending Stein's
method to
new distributions and applications.

There have been several previous developments of Stein's method using fixed
points of distributional transformations. \citet{Goldstein1997}
develop Stein's
method using the zero-bias transformation for which the normal
distribution is a
fixed point. Letting $U$ be a uniform $(0,1)$ random variable independent
of all
else, \citet{GoldsteinPC} and \citet{Pekoz2011} develop Stein's method
for the
exponential distribution using the fact that $W$ and $U W'$ have the same
distribution if and only if $W$ has an exponential distribution
[\citet{Pakes1992} and \citet{Lyons1995} also use this property]. We
will show below
that $W$ and $U W''$ have the same distribution if and only if $W$ is
distributed as the absolute value of a standard normal random variable
(see also
Remark~\ref{rem1} above). In this
light this paper can be viewed as extending the use of these types of
distributional transformations in Stein's method.

The layout of the remainder of the article is as follows. In
Section~\ref{sec1}
we discuss various properties and alternative representations of $K_s$, in
Section~\ref{sec2} we develop Stein's method for $K_s$ and prove
Theorem~\ref{thm2} and in Section~\ref{sec3} we prove Theorem~\ref
{thm1} by
constructing the coupling needed to apply Theorem~\ref{thm2} and
bounding the
appropriate terms.

\section{The distribution \texorpdfstring{$K_s$}{Ks}}\label{sec1}

In this section we collect some facts about $K_s$. Recall the notation and
definitions associated to the formula \eq{1} for the density $\kappa_s(x)$.
From \citeauthor{Abramowitz1992} [(\citeyear{Abramowitz1992}), Chapter 13], the Kummer $U$ function,
denoted $U(a,b,z)$, is the unique solution of
the differential equation
\[
z \frac{d^2 U}{dz^2}+(b-z) \frac{dU}{dz}-a U=0,
\]
which satisfies \eq{11} below. The
following lemma collects some facts about $U(a,b,z)$;
the {right} italic labeling of the formulas corresponds
to the equation numbers from \citeauthor{Abramowitz1992} [(\citeyear{Abramowitz1992}), Chapter 13],
and the notation
$U'(a,b,z)$ refers to the derivative with
respect to $z$.

%
\begin{lemma}\label{lem1} Let $a,b,z\in\IR$;
%
%
\begin{eqnarray}\qquad
&&\hspace*{6pt}\mbox{if $z>0$,}\qquad
U(a,b,z)=z^{1-b}U(1+a-b,2-b,z), \hspace*{36pt}\quad\textit{\aseq{1.29}}\label{5}
\\
\label{6}&&\hspace*{6pt}\mbox{if $a, z>0$,}
\nonumber
\\[-8pt]
&&\hspace*{6pt}\hspace*{293pt}\textit{\aseq{2.5}}
\\[-8pt]
\nonumber
&&\hspace*{6pt}\qquad U(a,b,z)=\frac{1}{\Gamma
(a)} \int
_0^\infty e^{-zt}t^{a-1}(1+t)^{b-a-1}
\,dt,
\\
&&\hspace*{6pt} U'(a,b,z) = -aU(a+1,b+1,z),\hspace*{104pt}\qquad\textit{\aseq{4.21}}\label{7}
\\
&&\hspace*{6pt}(1+a-b)U(a,b-1,z)
\nonumber
\\[-8pt]
&&\hspace*{6pt}\hspace*{287pt}\textit{\aseq{4.24}}
\\[-8pt]
\nonumber
&&\hspace*{6pt}\qquad=(1-b)U(a,b,z)-zU'(a,b,z), \label{8}
\\
&&\hspace*{6pt} U(a,b,z)-U'(a,b,z)=U(a,b+1,z), \hspace*{91pt}\quad\textit{\aseq{4.25}}\label{9}
\\
&&\hspace*{6pt}\qquad U(a-1,b-1,z) = (1-b+z)U(a,b,z)-zU'(a,b,z),\hspace*{6.5pt}\textit{\aseq{4.27}}\label{10}
\\
&&\hspace*{6pt} U(a,b,z) \sim z^{-a},\qquad(z\toinf), \hspace*{112pt}\quad\qquad\textit{\aseq{5.2}}\label{11}
\\
&&\hspace*{6pt}\mbox{for $a>-{ \frac{1}{2}}$}, \qquad U
\bigl(a,{ \tfrac{1}{2}},0 \bigr)=\Gamma\bigl({ \tfrac{1}{2}}
\bigr)/\Gamma\bigl(a+{ \tfrac{1}{2}} \bigr). \hspace*{68pt}\quad\textit{\aseq{5.10}}\label{12}
\end{eqnarray}
As a direct consequence of \eq{9} we have
%
%
\begin{equation}
\label{13} { \frac{\partial}{\partial z}} \bigl(e^{-z}U(a,b,z) \bigr) =
-e^{-z}U(a,b+1,z),
\end{equation}
combining \eq{7} and \eq{11} with $a=0$ we find
%
%
\begin{equation}
\label{14} U(0,b,z) = 1,
\end{equation}
and using \eq{5} with $a=-1/2$, $b=1/2$ and \eq{14} implies that for $z>0$,
%
%
\begin{equation}
\label{15} U \bigl(-{ \tfrac{1}{2}},{ \tfrac{1}{2}},z^2
\bigr) = z.
\end{equation}
\end{lemma}

By comparing integrands in \eq{6}, we also find the following fact.
%
%
\begin{lemma}\label{lem2}
Let $0<a<a'$, $b<b'$ and $z>0$. Then
%
%
\[
\qquad\Gamma(a) U(a, b, z) > \Gamma\bigl(a'\bigr) U\bigl(a',
b, z\bigr) \quad\mbox{and}\quad U(a,b,z)<U\bigl(a, b', z\bigr).
\]
\end{lemma}

The next results provide simpler representations for $K_s$.
%
%
\begin{proposition}\label{prop1} If $X$ and $Y$ are two independent random
variables having distributions
%
%
\[
X\sim\cases{ B(1,s-1), &\quad$\mbox{if $s>1$,}$\vspace*{2pt}
\cr
B(1/2,s-1/2),&\quad$\mbox{if
$1/2<s\leq1$},$}
\]
where $B(a,b)$ denotes the beta distribution,
and
%
%
\[
Y\sim\cases{ \Gamma(1/2,1 ), &\quad$\mbox{if $s>1$,}$\vspace*{2pt}
\cr
\Exp(1),&\quad
$\mbox{if $1/2<s\leq1$},$}
\]
where $\Gamma(a,b)$ denotes the gamma distribution and\/ $\Exp
(\lambda)$ the
exponential distribution, then
%
%
\[
\sqrt{2sXY}\sim K_s.
\]
\end{proposition}
\begin{pf}
Let $s>1$ and observe that by
first conditioning
on $X$,
we can express the density of $\sqrt{2sXY}$ as
%
%
\begin{equation}
\label{16} p_s(x):=\frac{\sqrt{2}(s-1)}{\sqrt{s\pi}} \int_0^1
\exp\biggl({ \frac{-x^2}{2sy}} \biggr) y^{-1/2}(1-y)^{s-2}
\,dy.
\end{equation}
After making the change of variable $y=1/(1+t)$ in \eq{16}, we find
%
%
\[
p_s(x)=\frac{\sqrt{2}(s-1)}{\sqrt{s\pi}}\int_0^\infty
\exp\biggl({ \frac{-x^2(t+1)}{2s}} \biggr)t^{s-2}(1+t)^{1/2-s}
\,dt,
\]
and now using \eq{6} with $a=s-1$ and $b=1/2$ in the definition \eq
{1} of
$\kappa_s$ implies that
$\kappa_s=p_s$.

Similarly, if $1/2<s\leq1$, then we can express the density of $\sqrt
{2sXY}$ as
%
%
\begin{equation}
\label{17} q_s(x):=\frac{ \Gamma(s) x}{s\sqrt{\pi}\Gamma(s-{ {1}/{2}})}\int_0^1
\exp\biggl({ \frac{-x^2}{2sy}} \biggr) y^{-3/2}(1-y)^{s-3/2}
\,dy,
\end{equation}
and after making the change of variable $y=1/(1+t)$ in \eq{17}, we find
\begin{eqnarray*}
q_s(x)&=&\frac{ \Gamma(s) x}{s\sqrt{\pi}\Gamma(s-{ {1}/{2}})}\int_0^\infty
\exp\biggl({ \frac{-x^2(t+1)}{2s}} \biggr)t^{s-3/2}(1+t)^{1-s}
\,dt
\\
&=&\Gamma(s)\sqrt{\frac{2}{s\pi}}\exp\biggl({ \frac
{-x^2}{2s}} \biggr)
\frac{x}{\sqrt{2s}}U \biggl(s-{ \frac{1}{2}}, { \frac{3}{2}}, {
\frac{x^2}{2s}} \biggr),
\end{eqnarray*}
where we have used \eq{6} with $a=s-1/2$ and $b=3/2$ in the second equality.
Applying \eq{5} with $a=s-1$ and $b=1/2$ to this
last expression implies $\kappa_s=q_s$.
\end{pf}

The previous representations easily yield useful formulas for Mellin transforms.
%
%
\begin{proposition}\label{prop2}
If $Z\sim K_s$ with $s\geq1/2$, then for all $r>-1$,
%
%
\begin{equation}
\label{18} \IE Z^r = \biggl(\frac{s}{2} \biggr)^{r/2}
\frac{\Gamma(s)\Gamma(r+1)}{\Gamma({ {r}/{2}}+s )}.\vadjust{\goodbreak}
\end{equation}
\end{proposition}
\begin{pf}
For $s>1/2$, we use Proposition~\ref{prop1} and well-known formulas
for the
Mellin transforms
of the beta and gamma distributions to find
%
%
\begin{equation}
\label{19} \IE Z^r = (2s)^{r/2}\frac{\Gamma(s)\Gamma({ {r}/{2}}+1 )
\Gamma({ {r}/{2}}
+{ {1}/{2}} )}{\Gamma({
{r}/{2}}+s )\Gamma({ {1}/{2}} )}.
\end{equation}
An application of the gamma duplication formula yields
%
%
\[
\Gamma\biggl({ \frac{r}{2}}+1 \biggr)\Gamma\biggl({
\frac{r}{2}}+{ \frac{1}{2}} \biggr) =\Gamma\biggl({
\frac{1}{2}} \biggr)2^{-r}\Gamma(r+1),
\]
which combined with \eq{19} implies \eq{18} for the case $s>1/2$.

The case $s=1/2$ follows from Proposition~\ref{prop3}(i) below which
implies that if $\law(Y)=\Exp(1)$, then $Z\eqlaw\sqrt{Y}$. Now \eq{19}
easily follows from well-known Mellin transform formulas and thus
\eq{18} also follows.
\end{pf}

In a few special cases we can simplify and extend Proposition~\ref
{prop1}. Below
$K_s(x)$ denotes the distribution function of $K_s$.

%
\begin{proposition}\label{prop3} We have the following special cases
of $K_s$:
\begin{eqnarray*}
\textup{(i)} &&\quad\kappa_{1/2}(x) = 2xe^{-x^2},
\\
\textup{(ii)} &&\quad\kappa_{1}(x) = (2/\pi)^{1/2}e^{-x^2/2},
\\
\textup{(iii)} &&\quad\lim_{s\toinf}K_{s}(x) = 1-e^{-\sqrt{2}x}.
\end{eqnarray*}
\end{proposition}
\begin{pf} The identities (i) and (ii) are immediate from \eq{15} and
\eq{14}, respectively. Using Stirling's formula for the
gamma function to take the limit as $s\toinf$ for fixed $r$ in \eq{18}
yields the moments of $\Exp(\sqrt2)$ which proves (iii).
\end{pf}

%
\begin{remark}
As discussed below, the preferential attachment model we study is a
special case
of a generalized P\'olya triangular urn scheme as studied by \citet
{Janson2006}.
The limiting distributions in his Theorem 1.3$(v)$ with $\alpha=2$ and
$\delta=\gamma=1$ include $K_s$. In fact, \citet{Janson2006},
Example 3.1,
discusses these limits, but, with the exception of the case $s=1$, it
does not
appear that the decomposition of Proposition~\ref{prop1} has
previously been
exposed. On the other hand, up to a scaling factor, the moment formula of
\citet{Janson2006}, Theorem 1.7, simplifies to that of
Proposition~\ref{prop2}
for $K_s$. The distribution
$K_s$ also appears in this urn context in
Section 9 of the survey article \citet{Janson2010}.

Additionally, if $Z\sim K_s$, then $Z^2/(2s)\sim D(1,1/2;s)$ for
$s\geq1/2$,
where $D(a,b;c)$ is a \textit{Dufresne law} as defined in \citet{Chamayou1999}.
Dufresne laws are essentially a generalization of products of
independent beta
and gamma random variables.
\end{remark}

We now collect one more fact about $K_s$, which will also prove useful in
developing the Stein's method framework below.

%
\begin{lemma}[(Mills ratio for $K_s$)]\label{lem3} For every $x\geq0$
and $s\geq
1$,
%
%
\[
\frac{1}{\kappa_s(x)}\int_x^\infty
\kappa_s(y) \,dy \leq\min\biggl\{\sqrt{\frac{\pi}{2}},
\frac{s}{x} \biggr\} .
\]
\end{lemma}
\begin{pf}
Using the definition \eq{1} of $\kappa_s$,
making the change of variable ${ \frac{y^2}{2s}} = z$ and
then applying
\eq{5}
with $a=s-1$ and $b=1/2$,
\eq{13} with $a=s-1/2$ and $b=1/2$ and then \eq{11} with $a=s-1/2$,
we find
\begin{eqnarray*}
\int_x^\infty\kappa_s(y) \,dy&=&
\frac{\Gamma(s)}{\sqrt{\pi}} \int_{{ {x^2}/{(2s)}}}^\infty z^{-1/2}
\exp(-z)U \biggl(s-1,{ \frac{1}{2}} ,z \biggr) \,dz
\\
&=&\frac{\Gamma(s)}{\sqrt{\pi}}\int_{{
{x^2}/{(2s)}}}^\infty\exp(-z)U
\biggl(s-{ \frac{1}{2}},{ \frac
{3}{2}},z \biggr) \,dz
\\
&=&\frac{\Gamma(s)}{\sqrt{\pi}}\exp\biggl({ \frac
{-x^2}{2s}} \biggr) U \biggl(s-{
\frac{1}{2}},{ \frac
{1}{2}},{ \frac{x^2}{2s}} \biggr),
\end{eqnarray*}
so that
%
%
\begin{equation}
\label{20} \frac{1}{\kappa_s(x)}\int_x^\infty
\kappa_s(y) \,dy =\sqrt{\frac{s}{2}} \frac{U (s-{ {1}/{2}},{
{1}/{2}},{ {x^2}/{(2s)}} )} {
U (s-1,{ {1}/{2}},{ {x^2}/{(2s)}} )}.
\end{equation}
First note that by applying \eq{5} with $a=s-1$ and $b=1/2$ in the denominator
of the final expression of
\eq{20} we have
%
%
\begin{equation}
\label{21} \frac{1}{\kappa_s(x)}\int_x^\infty
\kappa_s(y) \,dy = \frac{s}{x} \frac{U (s-{ {1}/{2}},{ {1}/{2}},{
{x^2}/{(2s)}} )} {
U (s-{ {1}/{2}},{ {3}/{2}},{ {x^2}/{(2s)}} )} \leq
\frac{s}{x},
\end{equation}
where the inequality follows by Lemma~\ref{lem2}.

Now applying \eq{5} to \eq{20} both in the denominator as before
and in the numerator
with $a=s-1/2$ and $b=1/2$, we find
\[
\frac{1}{\kappa_s(x)}\int_x^\infty
\kappa_s(y) \,dy = \sqrt{\frac{s}{2}} \frac{U (s,{ {3}/{2}},{
{x^2}/{(2s)}} )} {
U (s-{ {1}/{2}},{ {3}/{2}},{ {x^2}/{(2s)}} )} \leq
\sqrt{\frac{s}{2}}\frac{\Gamma(s-{ {1}/{2}} )}{\Gamma(s)},
\]
where again the inequality follows by Lemma~\ref{lem2}. Now applying
Lemma~\ref{lem4} below to this last expression and combining with \eq
{21} yields
the lemma.~%
\end{pf}

%
\begin{lemma}\label{lem4}
If $s\geq1$, then
\[
1 < \frac{\sqrt{s}\Gamma(s-{ {1}/{2}}
)}{\Gamma(s)}\leq\sqrt{\pi}.
\]
\end{lemma}
\begin{pf}
\citet{Bustoz1986}, Theorem 1, implies that
%
%
\begin{equation}
\label{22} \frac{\sqrt{s}\Gamma(s-{ {1}/{2}}
)}{\Gamma(s)}\vadjust{\goodbreak}
\end{equation}
is a decreasing function on $(1/2,\infty)$, so that for $s\geq1$, \eq
{22} is
bounded above by $\sqrt{\pi}$. Moreover, Stirling's formula implies
\[
\lim_{s\to\infty} \frac{ \Gamma(s)}{\sqrt{s} \Gamma
(s-{ {1}/{2}} )}=1.
\]
\upqed\end{pf}

\section{Stein's method for \texorpdfstring{$K_s$}{Ks}}\label{sec2}

In this section we develop Stein's method for $K_s$ and prove
Theorem~\ref{thm2}.

%
\begin{lemma}[(Characterizing Stein operator)]\label{lem5} If $Z\sim
K_s$ for
$s\geq1/2$, then for every twice differentiable function $f$ with
$f(0)=f'(0)=0$ and such that $\IE\vert f''(Z)\vert$, $\IE\vert
Zf'(Z)\vert$ and
$\IE\vert f(Z)\vert$ are finite, we have
%
%
\begin{equation}
\IE\bigl\{sf''(Z)-Zf'(Z)-2(s-1)f(Z)
\bigr\} = 0.\label{23}
\end{equation}
\end{lemma}
\begin{pf}
Let $C_s:=\sqrt{2}\Gamma(s)/\sqrt{s \pi}$. First note that
%
%
\begin{equation}
\label{24} \IE\bigl\{sf''(Z) \bigr\} =
C_s\int_{0}^\infty sf''(x)
\exp\biggl({ \frac{-x^2}{2s}} \biggr)U \biggl(s-1,{ \frac{1}{2}},{
\frac{x^2}{2s}} \biggr)\,dx.
\end{equation}
Using \eq{11} and \eq{7} with $a=s-1$ and $b=1/2$
we find that \eq{24} equals
\begin{eqnarray*}
&&C_s\int_{0}^\infty\!
f''(x) \int_{x}^\infty\! t
\exp\biggl({ \frac{-t^2}{2s}} \biggr) \biggl(U \biggl(s-1,{
\frac{1}{2}},{ \frac
{t^2}{2s}} \biggr)+(s-1) U \biggl(s,{
\frac{3}{2}}, { \frac{t^2}{2s}} \biggr) \biggr) \,dt \,dx
\\
&&\qquad= C_s\int_{0}^\infty
f'(t) t\exp\biggl({ \frac
{-t^2}{2s}} \biggr) \biggl(U
\biggl(s-1,{ \frac{1}{2}},{ \frac
{t^2}{2s}} \biggr)+(s-1) U \biggl(s,{
\frac{3}{2}}, { \frac{t^2}{2s}} \biggr) \biggr) \,dt
\\
&&\qquad= \IE\bigl\{Zf'(Z) \bigr\} + C_s\int
_{0}^\infty f'(t)\cdot(s-1) t\exp
\biggl({ \frac{-t^2}{2s}} \biggr) U \biggl(s,{ \frac{3}{2}},{
\frac{t^2}{2s}} \biggr) \,dt,
\end{eqnarray*}
where in the first equality we have used Fubini's theorem [justified
by\break
$\IE\vert f''(Z)\vert<\infty$] and the fact that
$f'(0)=0$.

We also have
\begin{eqnarray*}
&&\hspace*{-4pt} C_s\int_{0}^\infty f'(t)
\cdot t\exp\biggl({ \frac
{-t^2}{2s}} \biggr) (s-1)U \biggl(s,{
\frac{3}{2}},{ \frac
{t^2}{2s}} \biggr) \,dt
\\
&&\hspace*{-7pt}\qquad=C_s\int_0^\infty \!f'(t)
\int_t^\infty2(s-1)\!\exp\biggl({
\frac{-x^2}{2s}} \biggr)
\\
&&\hspace*{69pt}\qquad\quad{}\times\biggl( \biggl(-{ \frac{1}{2}}+{ \frac
{x^2}{2s}} \biggr)U
\biggl(s,{ \frac{3}{2}}, { \frac{x^2}{2s}} \biggr) -{
\frac{x^2}{2s}}U' \biggl(s,{ \frac
{3}{2}},{
\frac{x^2}{2s}} \biggr) \biggr) \,dx \,dt
\\
&&\hspace*{-7pt}\qquad= C_s\int_{0}^\infty f(x)\cdot2(s-1)
\exp\biggl({ \frac{-x^2}{2s}} \biggr)
\\
&&\hspace*{24pt}\qquad\quad{}\times\biggl( \biggl(-{ \frac{1}{2}}+{ \frac
{x^2}{2s}} \biggr)U
\biggl(s,{ \frac{3}{2}}, { \frac{x^2}{2s}} \biggr) -{
\frac{x^2}{2s}}U' \biggl(s,{ \frac
{3}{2}},{
\frac{x^2}{2s}} \biggr) \biggr) \,dx
\\
&&\hspace*{-7pt}\qquad= C_s\int_{0}^\infty f(x)\cdot2(s-1)
\exp\biggl({ \frac
{-x^2}{2s}} \biggr) U \biggl(s-1,{ \frac{1}{2}},{
\frac{x^2}{2s}} \biggr) \,dx
\\
&&\hspace*{-7pt}\qquad= \IE\bigl\{2(s-1)f(Z) \bigr\} ,
\end{eqnarray*}
where in the second equality we have used Fubini's theorem [justified
by $\IE\vert Zf'(Z)\vert<\infty$]
and the fact that $f(0)=0$, and in the third we have
used \eq{10} with $a=s$ and $b=3/2$.
Hence,
%
%
\[
\IE\bigl\{sf''(Z) \bigr\} = \IE\bigl
\{Zf'(Z) \bigr\} + \IE\bigl\{ 2(s-1)f(Z) \bigr\} ,
\]
which proves the claim.
\end{pf}

For the sake of brevity, let $V_s(x):=U (s-1,{ \frac
{1}{2}},{ \frac{x^2}{2s}} )$.

%
\begin{lemma}\label{lem6} For all functions $h$ such that $\IE h(Z)$ exists,
the second order differential equation
%
%
\begin{equation}
\label{25} sf''(x) - xf'(x) -
2(s-1)f(x) = h(x)-\IE h(Z)
\end{equation}
with initial conditions $f(0)=f'(0)=0$
has solution
%
%
\begin{eqnarray}
\label{26}
f(x)& =& \frac{1}{s} V_s(x) \int
_{0}^x\frac{1}{V_s(y)\kappa_s(y)}\int_{0}^y
\tilde{h}(z)\kappa_s(z)\,dz\,dy
\nonumber
\\[-8pt]
\\[-8pt]
\nonumber
&=& -\frac{1}{s}V_s(x) \int_{0}^x
\frac{1}{V_s(y) \kappa_s(y)}\int_{y}^\infty\tilde{h}(z)
\kappa_s(z)\,dz\,dy,
\end{eqnarray}
where $\tilde{h} = h - \IE h(Z)$.
\end{lemma}
In order to prove Lemma~\ref{lem6}, we use the following intermediate result.
%
%
\begin{lemma}\label{lem7}
If $g$ and $f$ are functions such that
$g(0)=f(0)=0$ and for $x>0$,
%
%
\begin{equation}
\label{27} sg'(x) - s\biggl({ \frac{x}{s}}-d(x)\biggr) g(x)
= \tilde{h}(x),\qquad f'(x) - d(x) f(x) = g(x),
\end{equation}
where
%
%
\begin{equation}
\label{28} d(x) = \frac{\partial}{\partial x} \log V_s(x) =
\frac{V_s'(x)}{V_s(x)},
\end{equation}
then $f$ solves \eq{25} and $f'(0)=0$.

Conversely, if $f$ is a solution to \eq{25} with $f(0)=f'(0)=0$
and $g(x)=f'(x)-d(x)f(x)$, then $g(0)=0$ and $f$ and $g$ satisfy $\eq{27}$.
\end{lemma}
\begin{pf}
Assume $f$ and $g$ satisfy \eq{27} and $f(0)=g(0)=0$. The fact that
$f'(0)=0$ follows easily from the second equation
of \eq{27}.
To show that \eq{27} yields a solution to \eq{25},
differentiate the second equality in \eq{27} and
combine the resulting equations to obtain
%
%
\[
sf''(x)-xf'(x)- \bigl(sd'(x)+sd(x)^2-xd(x)
\bigr) f(x) = \tilde{h}(x).\vadjust{\goodbreak}
\]
Hence, we only need to show that
%
%
\begin{equation}
\label{29} sd'(x)+sd(x)^2-xd(x) = 2(s-1).
\end{equation}
In order to simplify the calculations, let us introduce
\[
D(z) = \frac{\partial}{\partial z}\log U \biggl(s-1,{ \frac{1}{2}},z \biggr
) =
\frac{U' (s-1,{ {1}/{2}},z )}{U (s-1,
{ {1}/{2}},z )};
\]
note that $d(x)={ \frac{x}{s}}D ({ \frac
{x^2}{2s}} )$. With this and
$z={ \frac{x^2}{2s}}$, \eq{29} becomes
%
%
\begin{equation}
\label{30} \bigl({ \tfrac{1}{2}}-z \bigr) D(z) + zD'(z) +
zD(z)^2 = s-1.
\end{equation}
The left-hand side of \eq{30} is equal to
\begin{eqnarray*}
&&\frac{ ({ {1}/{2}}-z )U'
(s-1,{ {1}/{2}},z )
+zU'' (s-1,{ {1}/{2}},z )}{U
(s-1,{ {1}/{2}},z ) }
\\
&&\qquad= (s-1)\frac{ (-{ {1}/{2}}+z )U
(s,{ {3}/{2}},z )
-zU' (s,{ {3}/{2}},z )} { U (s-1,
{ {1}/{2}},z )} = s-1,
\end{eqnarray*}
where we have used \eq{7} with $a=s-1$ and $b=1/2$ to handle the
derivatives in the first equality
and then \eq{10} with $a=s$
and $b=3/2$ in the second. Hence, \eq{29} holds, as desired.

If $f$ is a solution to \eq{25} with $f(0)=f'(0)=0$
and $g(x)=f'(x)-d(x)f(x)$, then obviously $g(0)=0$
and the second assertion of the lemma follows from the
previous calculations.
\end{pf}
\begin{pf*}{Proof of Lemma~\ref{lem6}}
Lemma~\ref{lem7} implies that we
only need to solve~\eq{27}. Note first that the general differential
equation
%
%
\[
F'(x) - A'(x)F(x) = H(x), \qquad x > 0, F(0)=0,
\]
has solution
%
%
\[
F(x) = e^{A(x)}\int_{0}^x
H(z)e^{-A(z)}\,dz.
\]
Hence, noticing that
%
%
\[
\frac{x}{s}-d(x) = -\frac{\partial}{\partial x} \log\kappa_s(x),
\]
the solution to the first equation in \eq{27} is
%
%
\begin{equation}
g(y) = \frac{1}{\kappa_s(y)}\int_0^y
\frac{\tilde{h}(z)}{s} \kappa_s(z)\,dz, \label{31}
\end{equation}
whereas the solution to the second equation in \eq{27} is
%
%
\[
f(x) = V_s(x)\int_0^x
\frac{g(y)}{V_s(y)}\,dy,
\]
which is the first identity of \eqref{26}; the second follows by
observing that
$\int_0^\infty\tilde{h}(x)\kappa_s(x) \,dx = 0$.\vadjust{\goodbreak}
\end{pf*}

Before developing the Stein's method machinery
further we need two more lemmas, the first of which
is well known and easily read from \citet{Gordon1941}.

%
\begin{lemma}[(Gaussian Mills ratio)]\label{lem8}
For $x,s >0$,
\[
\exp\biggl({ \frac{x^2}{2s}} \biggr)\int_x^\infty
\exp\biggl({ \frac{-t^2}{2s}} \biggr)\,dt \leq\min\biggl\{\sqrt{
\frac{s\pi}{2}}, \frac{s}{x} \biggr\}.
\]
\end{lemma}

%
\begin{lemma}\label{lem9}
If $d(x)$ is defined by \eq{28}, then for $s\geq1$ and $x>0$
\begin{eqnarray*}
0&\leq&-d(x)\leq\frac{\sqrt{2} \Gamma(s)}{\sqrt{s} \Gamma
(s-{ {1}/{2}} )} < \sqrt{2},
\\
0&\leq&-x d(x)\leq2(s-1).
\end{eqnarray*}
\end{lemma}
\begin{pf}
To prove the first assertion note that \eq{7} with $a=s-1$ and $b=1/2$
followed by \eq{5}
with $a=s-1/2$ and $b=1/2$ and Lemma~\ref{lem2} imply
%
%
\begin{eqnarray}\label{32}
-d(x)&=&-\frac{x}{s}\frac{U' (s-1,{ {1}/{2}},
{ {x^2}/{(2s)}} )} { U (s-1,{ {1}/{2}},{ {x^2}/{(2s)}} ) } \nonumber\\
&=&\frac{\sqrt
{2}(s-1)}{\sqrt{s}}
\frac{U (s-{ {1}/{2}},{ {1}/{2}},
{ {x^2}/{(2s)}} ) }{ U (s-1,{ {1}/{2}},{ {x^2}/{(2s)}} ) }
\\
&\leq&\frac{\sqrt{2}(s-1)\Gamma(s-1)}{\sqrt{s} \Gamma
(s-{ {1}/{2}} )}.\nonumber
\end{eqnarray}
The claimed upper bound now follows from Lemma~\ref{lem4}.
The lower bound follows from the final expression of \eq{32},
since for $s>1$, the integral representation \eq{6} implies all terms
in the
quotient
are nonnegative, and for $s=1$, \eq{14} implies $d(x)=0$.

For the second assertion, we use \eq{8} with $a=s-1$ and $b=1/2$ in
the second
equality below to find
%
%
\begin{eqnarray}
\label{33} -xd(x) &=&-\frac{x^2}{s}\frac{U' (s-1,{ {1}/{2}},{
{x^2}/{(2s)}} )} {
U (s-1,{ {1}/{2}},{ {x^2}/{(2s)}}
) }
\nonumber
\\[-8pt]
\\[-8pt]
\nonumber
&=&2\biggl(s-{
\frac{1}{2}}\biggr)\frac{U (s-1,-{
{1}/{2}},{ {x^2}/{(2s)}} ) }{
U (s-1,{ {1}/{2}},{ {x^2}/{(2s)}}
) } -1.
\end{eqnarray}
Applying Lemma~\ref{lem2} to \eq{33} proves the remaining upper bound.
The second lower bound follows from the first.
\end{pf}

%
\begin{lemma}\label{lem10}
If $g$ satisfies the first equation of \eq{27} with $g(0)=0$, then
%
%
\[
g(x)=\frac{1}{s \kappa_s(x)}\int_0^x
\kappa_s(y) \tilde{h}(y) \,dy =-\frac{1}{ s \kappa_s(x)}\int
_x^\infty\kappa_s(y) \tilde{h}(y) \,dy.
\]
\begin{itemize}
\item If $h$ is nonnegative and bounded, then for all $x>0$ and $s\geq1$,
%
%
\begin{equation}
\label{34} \bigl\vert g(x)\bigr\vert\leq\Vert h\Vert\min\biggl\{\frac{1}{s}
\sqrt{\frac{\pi}{2}}, \frac{1}{x} \biggr\}.
\end{equation}
\item If $h$ is absolutely continuous with bounded derivative, then for all
$s\geq1$
%
%
\begin{equation}
\label{35} \Vert g\Vert\leq\bigl\Vert h'\bigr\Vert\biggl(1+
\frac{1}{s}\sqrt{\frac
{\pi
}{2}} \biggr).
\end{equation}
\end{itemize}
\end{lemma}
\begin{pf}
The first assertion is a restatement of \eq{31}, recorded in this
lemma for
convenient future
reference.

If $h(x)\geq0$ for $x\geq0$ with $\Vert h\Vert<\infty$, then for all
$s\geq1$ and
$x>0$,
%
%
\[
\bigl\vert g(x)\bigr\vert\leq\frac{\Vert\tilde{h}\Vert}{ s \kappa
_s(x)}\int_x^\infty
\kappa_s(y) \,dy \leq\min\biggl\{\sqrt{\frac{\pi}{2}},
\frac{s}{x} \biggr\} \frac
{\Vert h\Vert}{s},
\]
where we have used Lemma~\ref{lem3}; this shows \eq{34}.

Let $h$ be absolutely continuous with $\Vert h'\Vert<\infty$,
and without loss of generality assume that
$h(0)=0$ so that for $x\geq0$, $\vert h(x)\vert\leq\Vert h'\Vert
x$. In particular,
if $Z_s\sim K_s$, then
$\tilde{h}(x)\leq(x+\IE Z_s) \Vert h'\Vert$ and noting
that $\IE Z_s\leq\sqrt{\IE Z_s^2 }
=1$ (using Proposition~\ref{prop2}), we can apply Lemma~\ref{lem3} to find
that for $x>0$,
%
%
\[
\bigl\vert g(x)\bigr\vert\leq\frac{\Vert h'\Vert}{s \kappa_s(x)} \int
_x^\infty(y+1)
\kappa_s(y) \,dy \leq\frac{\Vert h'\Vert}{s} \biggl(\frac{\int
_x^\infty y\kappa_s(y)\,dy}{\kappa_s(x)} +
\sqrt{\frac{\pi}{2}} \biggr).
\]
To bound the integral in this last expression, we
make the change of variable ${ \frac{y^2}{2s}} = z$ and
apply \eq{13} with
$a=s-1$, $b=-1/2$ and \eq{11} with $a=s-1$ to find
%
%
\begin{eqnarray*}
\frac{\int_x^\infty y \kappa_s(y) \,dy}{\kappa_s(x)}& =&\frac{s}{\kappa
_s(x)} \int_{ {x^2}/{(2s)}}^\infty
e^{-z} U \biggl(s-1, { \frac{1}{2}}, z \biggr) \,dz\\
& =&s
\frac{U (s-1, -{ {1}/{2}}, { {x^2}/{(2s)}} )} {
U (s-1, { {1}/{2}}, { {x^2}/{(2s)}} )} \leq s,
\end{eqnarray*}
where the last inequality follows from Lemma~\ref{lem2}.
\end{pf}

%
\begin{lemma}\label{lem11}
Let $f$ be defined as in~\eq{26} with $f(0)=f'(0)=0$.
\begin{itemize}
\item If $h$ is nonnegative and bounded and $s\geq1$, then
%
%
\begin{equation}
\label{36} \bigl\Vert f'\bigr\Vert\leq\sqrt{2\pi}\Vert h\Vert.
\end{equation}
\item If $h$ is nonnegative, bounded and absolutely continuous with bounded
derivative and $s\geq1$, then
%
%
\begin{equation}
\bigl\Vert f''\bigr\Vert\leq2 \biggl(\pi\sqrt{s}+
\frac{1}{s} \biggr)\Vert h\Vert. \label{37}
\end{equation}
If $s=1/2$, then
%
%
\begin{equation}
\label{38} \bigl\Vert f''\bigr\Vert\leq4\Vert h\Vert.\vadjust{\goodbreak}
\end{equation}
\item If $h$ is absolutely continuous with bounded derivative and
$s\geq1$,
then
%
%
\begin{equation}
\label{39} \bigl\Vert f'''\bigr\Vert\leq8
\biggl(s+\frac{1}{4}+\sqrt{\frac{\pi
}{2}} \biggr)\bigl\Vert h'
\bigr\Vert.
\end{equation}
If $s=1/2$, then
%
%
\begin{equation}
\label{40} \bigl\Vert f'''\bigr\Vert\leq4\bigl\Vert
h'\bigr\Vert.
\end{equation}
\end{itemize}
\end{lemma}
\begin{pf}
From \eq{26} of Lemma~\ref{lem6} we have that
%
%
\[
f(x)=V_s(x)\int_0^x
\frac{g(y)}{V_s(y)} \,dy,
\]
where $g$ is as in Lemma~\ref{lem10}.
If either $h$ is bounded or absolutely continuous with bounded
derivative, then
recall that Lemma~\ref{lem10} implies $g$ is bounded.
If $s\geq1$, then~\eq{7} and \eq{11} with $a=s-1$ and $b=1/2$ imply that
$V_s(x)=U(s-1,{ \frac{1}{2}}, { \frac{x^2}{2s}})$ is
nonincreasing and positive for positive $x$, so that
%
%
\begin{equation}
\label{41} \bigl\vert f(x)\bigr\vert\leq x\Vert g\Vert.
\end{equation}
Now, again by \eq{27}, we have
%
%
\begin{equation}
\label{42} \bigl\vert f'(x)\bigr\vert\leq\bigl\vert d(x)f(x)\bigr\vert
+\Vert g
\Vert\leq\Vert g\Vert\bigl(\bigl\vert xd(x)\bigr\vert+1 \bigr) \leq
\Vert g\Vert(2s-1),
\end{equation}
where we have used \eq{41} in the first inequality and Lemma~\ref
{lem9} in the
second. Applying the bound \eq{34} proves \eq{36}.

To bound $f''$ for $h$ having $\Vert h'\Vert<\infty$, let $s\geq
1/2$ and
differentiate \eq{25} to find
%
%
\begin{equation}
\label{43} f'''(x)-{
\frac{x}{s}}f''(x)=\frac{2s-1}{s}
f'(x) + \frac
{h'(x)}{s},
\end{equation}
which implies
%
%
\[
\frac{d}{dx} \biggl(\exp\biggl({ \frac{-x^2}{2s}}
\biggr)f''(x) \biggr) =\exp\biggl({
\frac{-x^2}{2s}} \biggr) \biggl(\frac
{2s-1}{s} f'(x) +
\frac{h'(x)}{s} \biggr).
\]
Integrating, we obtain
%
%
\[
\exp\biggl({ \frac{-x^2}{2s}} \biggr)f''(x) =-
\int_x^\infty\exp\biggl({ \frac{-y^2}{2s}}
\biggr) \biggl(\frac
{2s-1}{s} f'(y) + \frac{h'(y)}{s} \biggr)
\,dy,
\]
so that Lemma~\ref{lem8} yields
%
%
\begin{eqnarray}
\label{44} \bigl\vert f''(x)\bigr\vert&\leq&(2s-1)\bigl\Vert
f'\bigr\Vert\min\biggl\{\sqrt{\frac
{\pi
}{2s}}, \frac{1}{x}
\biggr\}
\nonumber
\\[-8pt]
\\[-8pt]
\nonumber
&&{} + \frac{1}{s}\exp\biggl({ \frac{x^2}{2s}} \biggr)\int
_x^\infty\exp\biggl({ \frac{-y^2}{2s}}
\biggr)h'(y)\,dy.
\end{eqnarray}
If $\Vert h\Vert<\infty$, then an integration by parts yields a
bound on
the second
term of~\eq{44} which yields
%
%
\[
\bigl\vert f''(x)\bigr\vert\leq(2s-1)\bigl\Vert f'
\bigr\Vert\min\biggl\{\sqrt{\frac{\pi}{2s}}, \frac{1}{x} \biggr\} +
\frac{2\Vert h\Vert}{s}.
\]
If $s\geq1$, then apply the bound \eq{36} above on $\Vert f'\Vert$ to
find \eq{37};
for $s=1/2$, \eq{38} follows immediately. Now, we can apply
Lemma~\ref{lem8} directly to \eq{44} to find
%
%
\begin{equation}
\label{45} \bigl\vert f''(x)\bigr\vert\leq\bigl((2s-1)\bigl\Vert
f'\bigr\Vert+\bigl\Vert h'\bigr\Vert\bigr) \min\biggl\{\sqrt{
\frac{\pi}{2s}}, \frac{1}{x} \biggr\}.
\end{equation}
Finally, \eq{43} implies
%
%
\begin{equation}
\label{46} s\bigl\vert f'''(x)\bigr\vert\leq
\bigl\vert xf''(x)\bigr\vert+(2s-1)\bigl\Vert f'\bigr\Vert
+\bigl\Vert h'\bigr\Vert;
\end{equation}
the first term can be bounded by \eq{45}, and if $s\geq1$, a subsequent
application of~\eq{42} on $\Vert f'\Vert$ and then \eq{35} on
$\Vert g\Vert$
yields~\eq{39}. If $s=1/2$, then \eq{40} follows from \eq{46} and
\eq{45}.
\end{pf}

In order to obtain the bounds for the Kolmogorov metric, we need to introduce
the smoothed half-line indicator function
%
%
\begin{equation}
\label{47} h_{a,\eps}(x) = \frac{1}{\eps}\int_0^\eps
\I[x\leq a+t] \,dt.
\end{equation}

%
\begin{lemma}\label{lem12}
If $Z\sim K_s$ and $W$ is a nonnegative random variable and $s\geq1$, then,
for all $\eps>0$,
%
%
\[
\dk\bigl(\law(W), K_s \bigr) \leq\sup_{a\geq0}\bigl\vert\IE
h_{a,\eps}(W)-\IE h_{a,\eps}(Z)\bigr\vert+\eps\sqrt{2}.
\]
If $s=1/2$, then, for all $\eps>0$,
%
%
\[
\dk\bigl(\law(W), K_{1/2} \bigr) \leq\sup_{a\geq0}\bigl\vert\IE
h_{a,\eps}(W)-\IE h_{a,\eps}(Z)\bigr\vert+\eps\sqrt{2/e}.
\]
\end{lemma}
\begin{pf}
The lemma follows from a well-known argument and the following bounds
on the
density $\kappa_s(x)$ defined by \eq{1}. If $s\geq1$, then by \eq
{7} with
$a=s-1$ and \eq{6}
with $a=s$, $\kappa_s(x)$ is nonincreasing in $x$ and
from \eq{12} with $a=s-1$,
%
%
\[
\kappa_s(0) =\frac{\Gamma(s)\sqrt{2}}{\Gamma(s-{{1}/{2}} )\sqrt{s}}
\leq\sqrt{2},
\]
where the inequality is by Lemma~\ref{lem4}. If $s=1/2$, then
$\kappa_s(x)=2xe^{-x^2}$ which has maximum $\sqrt{2/e}$.
\end{pf}

We will also need the following ``indirect'' concentration inequality; it
follows from the arguments of the proof of Lemma~\ref{lem12}
immediately above.

%
\begin{lemma}\label{lem13}
If $Z\sim K_s$ and $W$ is a nonnegative random variable and $s\geq1$, then,
for all\/ $0\leq a <b$,
%
%
\[
\IP(a < W \leq b)\leq\sqrt{2}(b-a) + 2\dk\bigl(\law(W), K_s
\bigr).
\]
If $s=1/2$, then, for all\/ $0\leq a <b$,
%
%
\[
\IP(a < W \leq b)\leq\sqrt{2/e}(b-a) + 2\dk\bigl(\law(W), K_s
\bigr).
\]
\end{lemma}

%
\begin{lemma}\label{lem14}
If $f$ satisfies \eq{25} for $h_{a,\eps}$ and $s\geq1$, then for
$x\geq0$,
\begin{eqnarray*}
s\bigl\vert f''(x+t)-f''(x)
\bigr\vert&\leq& \vert t\vert\bigl(2x(\pi\sqrt{s}+1)+(2s-1)\sqrt{2\pi}
\bigr)
\\
&&{}+\frac{1}{\eps}\int_{t \wedge0}^{t \vee0} \I[a< x+u\leq a+
\eps]\,du.
\end{eqnarray*}
If $s=1/2$, then for $x\geq0$,
\[
{ \frac{1}{2}}\bigl\vert f''(x+t)-f''(x)
\bigr\vert\leq4\vert t\vert x +\frac{1}{\eps}\int_{t \wedge0}^{t \vee0}
\I[a< x+u\leq a+\eps]\,du.
\]
\end{lemma}
\begin{pf}
Using \eq{25}, we obtain
\begin{eqnarray*}
s\bigl(f''(x+t)-f''(x)
\bigr) &=& x\bigl(f'(x+t)-f'(x)\bigr)+t
f'(x+t)
\\
&&{}+2(s-1) \bigl(f(x+t)-f(x)\bigr) + h_{a,\eps}(x+t)-h_{a,\eps}(x),
\end{eqnarray*}
hence,
\begin{eqnarray*}
s\bigl\vert f''(x+t)-f''(x)
\bigr\vert&\leq& \vert t\vert\bigl(x\bigl\Vert f''\bigr\Vert+\bigl
\Vert
f'\bigr\Vert+2(s-1)\bigl\Vert f'\bigr\Vert\bigr)
\\
&&{}+\frac{1}{\eps}\int_{t \wedge0}^{t \vee0} \I[a < x+u\leq
a+\eps]\,du.
\end{eqnarray*}
Applying the bounds of Lemma~\ref{lem11} yields the claim.
\end{pf}

%
\begin{lemma}\label{lem15}
Let $W$ be a nonnegative random variable with $\IE W^2 =1$ and let
$W^*$ be the
$s$-TDSB of\/ $W$ as in Definition~$\ref{def1}$ for some $s\geq1/2$.
For every
twice differentiable function $f$ with $f(0)=f'(0)=0$ and such that the
expectations below are well defined, we have
\[
s\IE f''\bigl(W^{*}\bigr)=\IE\bigl
\{Wf'(W)+2(s-1)f(W) \bigr\} .
\]
\end{lemma}
\begin{pf}
The lemma will follow from two facts:
\begin{itemize}
\item If $W''$ has the double size bias distribution of $W$, then for
all $g$
with $\IE\vert W^2 g(W)\vert<\infty$,
\[
\IE g\bigl(W''\bigr) =\IE\bigl\{W^2
g(W) \bigr\} .
\]
\item If $g$ is a function such that $g'(0)=g(0)=0$
and for $V$ as defined in Definition~\ref{def1}, $\IE\vert
g''(V)\vert
<\infty$, then
\[
s \IE g''(V)= g'(1)+2 (s-1) g(1).
\]
\end{itemize}
The first item above is easy to verify from the definition of the size bias
distribution and the fact that $\IE W^2=1$, and the second follows from
a simple
calculation after noting that $V$ has density
$ (2-{ \frac{1}{s}} )-2x (1-{ \frac
{1}{s}} )$ for $0<x<1$.

By conditioning on $W''$ and using the second fact above for
$g(t)=f(tW'')/\break  (W'')^2$, we find
%
%
\[
s\IE f''\bigl(W^{*}\bigr) = \IE\biggl\{
\frac{ f'(W'')}{W''}+2(s-1) \frac
{f(W'')}{(W'')^2} \biggr\} ,
\]
and applying the first fact above proves the lemma.
\end{pf}

\begin{pf*}{Proof of Wasserstein bound of Theorem~\ref{thm2}} Making
use of
Lem\-ma~\ref{lem6} and Lemma~\ref{lem15}, we have
\begin{eqnarray*}
\IE h(W) - \IE h(Z) &=& \IE\bigl\{sf''(W)-Wf'(W)-2(s-1)f(W)
\bigr\}
\\
&=& s\IE\bigl\{f''(W)-f''
\bigl(W^*\bigr) \bigr\} ,
\end{eqnarray*}
where $f$ is given by \eq{26}. If $h$ is Lipschitz continuous, then
$f$ is three
times differentiable almost everywhere and we have
\[
\bigl\vert\IE h(W) - \IE h(Z) \bigr\vert\leq s\bigl\Vert f'''
\bigr\Vert\IE\bigl\vert W-W^*\bigr\vert.
\]
We now obtain \eq{3} by invoking \eq{39} and \eq{40} of Lemma~\ref{lem11}.
\end{pf*}

\begin{pf*}{Proof of Kolmogorov bound of Theorem~\ref{thm2}} Fix
$a>0$ and let
$\eps>0$, to be chosen later. Let $f$ be as in \eq{26} with $\tilde
{h}$ replaced
by $h_{a,\eps}-\IE h_{a, \eps}(Z)$, where $h_{a, \eps}$ is defined
by \eq{47}. Define the indicator random variable $J =
\I[\vert W-W^*\vert\leq\beta]$. Now,
\begin{eqnarray*}
&&\IE h_{a,\eps}(W) - \IE h_{a,\eps}(Z)
\\
&&\qquad= s\IE\bigl\{f''(W)-f''
\bigl(W^*\bigr) \bigr\}
\\
&&\qquad= s\IE\bigl\{J\bigl(f''(W)-f''
\bigl(W^*\bigr)\bigr) \bigr\} + s\IE\bigl\{ (1-J) \bigl(f''(W)-f''
\bigl(W^*\bigr)\bigr) \bigr\}
\\
&&\qquad=: R_1 + R_2.
\end{eqnarray*}
If $s\geq1$, using \eq{37} from Lemma~\ref{lem11} implies
\[
\vert R_2\vert\leq4\bigl(\pi s^{3/2}+1\bigr)\IP\bigl(
\bigl\vert W-W^*\bigr\vert>\beta\bigr) \leq17s^{3/2} \IP\bigl(\bigl
\vert W-W^*\bigr\vert>
\beta\bigr).
\]
Applying Lemma~\ref{lem14},
\[
\vert R_1\vert\leq\beta\bigl(2\IE W(\pi\sqrt{s}+1)+(2s-1)\sqrt{2
\pi} \bigr) +\frac{1}{\eps}\int_{-\beta}^{\beta}\IP
(a< W+u\leq a+\eps)\,du.
\]
Noticing that $\IE W\leq1$ and applying Lemma~\ref{lem13} to the integrand,
\[
\vert R_1\vert\leq12s\beta+ 2\beta\eps^{-1}(\sqrt{2}
\eps+2\delta) \leq15s\beta+ 4\beta\eps^{-1}\delta,
\]
where $\delta= \dk(\law(W),K_s)$.

From Lemma~\ref{lem12}, we have
\[
\delta\leq\sqrt{2}\eps+ 15s\beta+ 4\beta\eps^{-1}\delta+
17s^{3/2} \IP\bigl(\bigl\vert W-W^*\bigr\vert>\beta\bigr).
\]
Choosing $\eps=8\beta$ and solving for $\delta$,
\[
\delta\leq16\sqrt{2}s\beta+ 30s\beta+ 34s^{3/2} \IP\bigl(\bigl\vert W-W^*
\bigr\vert>\beta\bigr),
\]
which yields \eq{4}.

A nearly identical argument yields the statement for $s=1/2$.
\end{pf*}

\section{\texorpdfstring{Proof of Theorem~\protect\ref{thm1}}{Proof of Theorem 1.1}}\label{sec3}

We first reformulate Theorem~\ref{thm1} in terms of a generalized P\'olya
urn model.
An urn initially contains $i$ black balls and $j$ white balls
and at each step a ball is drawn. If the
ball drawn is black, it is returned to the urn along with an additional
$\alpha$
black balls and $\beta$ white balls; if the ball drawn is white, the
ball is
returned to the urn along with an additional $\gamma$ black balls and
$\delta$
white
balls. We use the notation $(\alpha,\beta;\gamma,\delta)^n_{i,j}$ to
denote the distribution of the number of white balls in the urn after
$n$ draws
and replacements.
For example, $(\alpha,\beta;\gamma,\delta)^0_{i,j}$ has a single
point mass at
$j$ and
also note that $(1,0;0,1)_{i,j}^n$
corresponds to the classical P\'olya urn model.

\renewcommand{\thetheorema}{\protect\ref{thm1}\textup{(a)}}
%
\begin{theorema}\label{th1a}
Let $n\geq1$ and $i\geq0$ be integers and $\law(W_{n,i})= (2,0;\break1,1)^n_{i,1}$.
If $b_{n,i}^2=\IE W_{n,i}^2$, then, for some constants $c, C>0$
independent of~$n$,
\[
\frac{c}{\sqrt{n}}\leq\dk\bigl(\law(W_{n,i}/b_{n,i}),K_{(i+1)/2}
\bigr)\leq\frac{C}{\sqrt{n}}.
\]
\end{theorema}

Theorem~\ref{thm1} follows immediately from Theorem~\ref{thm1}(a) after noting
that
for model~1 with $n\geq i\geq2$,
the degree of vertex $i$ in $G_n$, the graph with $n$ vertices and
$n-1$ edges,
has distribution
$(2,0;1,1)^{n-i}_{2i-3,1}$; this is because the degree of vertex $i$ in
$G_i$ is
$1$ and
the sum
of the degrees of the remaining vertices is $2i-3$ (since $G_i$ has $i-1$
edges).
For model~2 with $n\geq i\geq1$, the degree of vertex $i$ in $G_n$,
the graph
with $n$ vertices and $n$ edges,
has distribution
$(2,0;1,1)^{n-i+1}_{2i-2,1}$; this is because the sum of the degrees of
$G_{i-1}$
is $2i-2$ and vertex $i$ has probability $1/(2i-1)$ of self-attachment when
forming
$G_i$ from $G_{i-1}$.

The lower bound of the theorem follows from the following general result
combined with the fact from Lemma \ref{lem21} below that $\IE
W_{n,i}^2\leq
2(1+2n)$.
%
%
\begin{lemma}
Let $\mu$ be a probability distribution with a density $f$ such that
for all $x$
in some interval $(a,b)$,
$f(x)>\eps>0$.
If $(X_n)_{n\geq1}$ is a sequence of integer-valued random variables and
$(a_n)_{n\geq1}$\vadjust{\goodbreak}
is a sequence of nonnegative numbers tending to zero, then
\[
\dk\bigl(\law(a_n X_n), \mu\bigr)\geq c
a_n
\]
for some positive constant $c$ independent of $n$.
\end{lemma}
\begin{pf}
Let $F$ be the distribution function of $\mu$ and note that the
hypothesis on the density $f$ implies that if $a\leq x <y \leq b$, then
%
%
\begin{equation}
F(y)-F(x)\geq\eps(y-x). \label{48}
\end{equation}
Since $\lim_n a_n=0$, there exists $N$ such that for all $n\geq N$,
there is an integer $k_n$ such that $[a_n k_n , a_n (k_n+1)]\subset(a,b)$.
From \eq{48}, for $n\geq N$ we have
\[
F\bigl(a_n (k_n+1)\bigr)-F(a_n
k_n)\geq a_n \eps,
\]
and now using the continuity of $F$ on $(a,b)$ and the fact that the
distribution
function $G_n$ of $a_n X_n$ is constant on the interval $I_n:=[a_n k_n
, a_n
(k_n+1))$,
it follows that for $n\geq N$,
\[
\dk\bigl(\law(a_n X_n), \mu\bigr)\geq
\sup_{x\in I_n}\bigl\vert G_n(x)-F(x)\bigr\vert\geq\frac{\eps}{2}
a_n.
\]
Since $G_n$ is the distribution function of a discrete random variable
and $F$
is continuous,
it follows that $\dk(\law(a_n X_n), \mu)>0$ for all
$n\in\IN$ (and in
particular $n<N$), so that we may choose $c>0$.
\end{pf}

%
\begin{remark}
As mentioned in the \hyperref[intr]{Introduction }we write our results
as rates, but the
constants are recoverable.
For the sake of clarity, we have not been careful to optimize
the bounds in our arguments, but it is clear that sharper statements
can be read from the proofs below.
For example, the constant in both the lower bound and upper bounds of
Theorem~\ref{thm1}(a)
depend crucially on the scaling factor $\IE W_{n,i}^2$. For our purposes
Lemma~\ref{lem21} below is acceptable, but
note that exact results are available [see \eq{70} and \eq{71} in the
proof of Lemma~\ref{lem21}].
\end{remark}

Now let $W:=W_{n,i}$ have distribution $(2,0;1,1)^n_{i,1}$.
We will use \eq{4} to prove the upper bounds of
Theorem~\ref{thm1}(a) and so we will show that
there is a close coupling of $W$ and
$VW''$, where $V$ is as in Definition~\ref{def1} with $s=(i+1)/2$.
This result
will follow from the following lemmas
proved at the end
of this section.

%
\begin{lemma}\label{lem16}
There is a coupling $(R,W'')$ of
$(2,0;1,1)^{n-1}_{i,3}$ and the double size bias distribution of
$(2,0;1,1)_{i,1}^n$
satisfying $\IP(R\neq W'')\leq C/\sqrt{n}.$
\end{lemma}

%
\begin{lemma}\label{lem17}
The distribution $(2,0;1,1)_{i,1}^{n}$ can be expressed as a mixture of the
distributions $(2,0;1,1)_{i+1,2}^{n-1}$ and
$(2,0;1,1)_{i+2,1}^{n-1}$ with respective probabilities $1/(1+i)$ and
$1-1/(1+i)$.
\end{lemma}

In the next lemma we use the notation $(\alpha,\beta;\gamma,\delta
)^N_{i,j}$ for a nonnegative integer-valued random variable $N$ to
denote a mixture of the distributions $(\alpha,\beta;\gamma,\delta
)^n_{i,j}$ for $n=0,1,2,\ldots$ that are mixed with respective
probabilities $\IP(N=n)$ for $n=0,1,2,\ldots.$
%
%
\begin{lemma}\label{lem18}
Let $\law(R)=(2,0;1,1)^{n-1}_{i,3}$, let $\law
(X_1)=(1,0;0,1)_{1,2}^{R-3}$ and let $\law
(X_2)=(1,0;0,1)_{2,1}^{R-3}$. Then
%
%
\[
\law(X_1)=(2,0;1,1)_{i+1,2}^{n-1} \quad\mbox{and}\quad
\law(X_2)= (2,0;1,1)_{i+2,1}^{n-1}.
\]
\end{lemma}

%
\begin{lemma}\label{lem19}
Let $U_1$ and $U_2$ be uniform $(0,1)$ random variables, independent of each
other and of $R$, defined
as in Lemma~$\ref{lem16}$. Then there exist random variables $X_1$ with
distribution
$(1,0;0,1)_{1,2}^{R-3}$ and $X_2$ with distribution
$(1,0;0,1)_{2,1}^{R-3}$ such
that
%
%
\[
\bigl\vert X_1 - R\max(U_1,U_2)\bigr\vert<3
\quad\mbox{and}\quad\bigl\vert X_2 - R\min(U_1,U_2)\bigr\vert<3\qquad
\mbox{a.s.}
\]
\end{lemma}

From these lemmas we can now prove Theorem~\ref{thm1}(a); here and below
we use $C$ to denote a generic constant that may differ from line to line.

\begin{pf*}{Proof of Theorem~\ref{th1a}}
Let $W=W_{n,i}$ and $b=b_{n,i}$,
let $(R,W'')$ be defined as in Lemma~\ref{lem16} above and, as per
Definition~\ref{def1},
let $Y$\break be~a~Bernoulli($1/(1+i)$) random variable
and $ V = Y \max(U_{1},U_{2})+\break(1-Y)\min(U_{1},U_{2})$, where
$U_1$ and $U_2$ are independent uniform $(0,1)$ variables independent
of $Y$.
Lemmas~\ref{lem17},~\ref{lem18} and~\ref{lem19} imply that
we can couple $W$ and $V R$ together so that $|W-VR|<3$ almost surely.
Thus, using Lemma~\ref{lem16},
%
%
\[
\IP\bigl(\bigl|W-VW''\bigr|>3\bigr)\leq\IP
\bigl(W''\not=R\bigr)\leq C/\sqrt{n},
\]
and recalling that $VW''$ has the $s$-TDSB distribution,
the theorem follows from~\eq{4}
taking $\beta= 3/b$, noting
that for $c>0$, $(cW)'\eqlaw cW'$ and using
$b^2\geq Cn$ from Lemma~\ref{lem21} below.
\end{pf*}

We have left to prove Lemmas~\ref{lem16}--\ref{lem19}
and~\ref{lem21}; Lemma~\ref{lem17} is immediate after considering the urn
process corresponding to $(2,0;1,1)^{n}_{i,1}$ and conditioning on the
color of the
first ball drawn, which is white with probability $1/(1+i)$.

\begin{pf*}{Proof of Lemma~\ref{lem18}}
Consider an urn with $i$ green balls, 1 black ball and 2 white
balls. A ball is drawn at random and replaced in the urn along with
another ball of
the same color plus an additional green ball.

If $X$ is the number
of times a nongreen ball is drawn in $n-1$ draws, the number of white
balls in
the urn
after $n-1$ draws is distributed as $(1,0;0,1)_{1,2}^{X}$.
Since $X+3$ is distributed as
$(2,0;1,1)_{i,3}^{n-1}$ (which is also that of $R$) and the number of white
balls in the urn after $n-1$ draws
has distribution $(2,0;1,1)_{i+1,2}^{n-1}$, the first equation follows. The
second equation
follows from similar considerations.
\end{pf*}

\begin{pf*}{Proof of Lemma~\ref{lem19}}
We will show that for $U_1$ and $U_2$ independent uniform $(0,1)$ random
variables,
there exist random variables $N$ and $M$ such that
$\law(N)=(1,0;0,1)_{1,2}^{n-3}$, $\law(M)=\law(n-N)$ and
\begin{eqnarray*}
\bigl|N-n\max(U_1,U_2)\bigr|<3 \quad\mbox{and}\quad\bigl|M-n
\min(U_1,U_2)\bigr|<3\qquad\mbox{a.s.}
\end{eqnarray*}
The lemma follows from these ``conditional'' almost sure statements
after noting that $\law(M)=(1,0;0,1)_{2,1}^{n-3}$ since we can think
of $n-N$ as the number of
black balls in the $(1,0;0,1)_{1,2}^{n-3}$ urn.

The formulas of \citeauthor{Durrett2010} [(\citeyear{Durrett2010}), page 206] imply that
$(1,0;0,1)_{1,2}^{n-3}$
has distribution
function
%
%
\begin{equation}
F(k) = \biggl(\frac{k}{n-1} \biggr) \biggl(\frac{k-1}{n-2} \biggr
),\qquad k=1,
\ldots,n-1, \label{49}
\end{equation}
and it is straightforward to verify
that
\[
N:= \max\bigl({\bigl\lceil(n-1)U_1\bigr\rceil}, 1+{\bigl
\lceil(n-2)U_2\bigr\rceil}\bigr)
\]
has the same distribution. We find $|N-n\max(U_1,U_2)|<3$ and thus a
coupling satisfying the first claim above.
Defining
\[
M:= \min\bigl({\bigl\lceil(n-1)U_1\bigr\rceil}, 1+{\bigl
\lceil(n-2)U_2\bigr\rceil}\bigr),
\]
\eq{49} implies $\law(M)=\law(n-N)$
and $|M-n\min(U_1,U_2)|<3$.
\end{pf*}

Before proving Lemma~\ref{lem16},
we provide a useful construction
for the double size bias distribution of a sum of indicators.

%
\begin{lemma}\label{lem20}
Let $W=\sum_{i=1}^n X_i$, where the $X_i$ are Bernoulli random
variables and
$b^2:=\IE W^2$. For each $j,k\in\{1,\ldots, n\}$, let
$ (X_i^{(j,k)} )_{i\notin\{j,k\}}$ have the distribution of
$(X_i)_{i\notin\{j,k\}}$
conditional on $X_j=X_k=1$ and let $J$ and $K$ be random variables independent
of the variables above satisfying
\[
\IP(J=j,K=k)=\frac{\IE(X_j X_k)}{b^2}, \qquad\mbox{$j,k\geq1$.}
\]
%
Then,
%
%
\[
W''=\sum_{i\notin\{J,K\}}
X_i^{(J,K)} +2-\I[J=K]
\]
has the double size bias distribution of $W$.
\end{lemma}
\begin{pf}
We have
\begin{eqnarray*}
\IE f\bigl(W''\bigr) &=&b^{-2}\sum
_{j,k}\IE(X_jX_k)\IE f \biggl(\sum
_{i\notin\{j,k\}} X_i^{(j,k)} +2-\I[j=k]
\biggr)
\\
&=&b^{-2}\sum_{j,k}\IE(X_jX_k)
\IE\bigl\{ f (W ) | X_j=X_k=1\bigr\}
\\
&=&b^{-2}\sum_{j,k}\IE\bigl
\{X_jX_k f (W )\bigr\}=b^{-2}\IE\bigl
\{W^2 f (W ) \bigr\} ;
\end{eqnarray*}
this is exactly \eq{2}.
\end{pf}

To simplify the notation we consider $i$ fixed in what follows. We write
\[
W_n=\sum_{j=0}^n
X_{j},
\]
where for $j\geq1$, $X_{j}$ is the indicator that a white ball is
drawn on
draw $j$ from the $(2,0;1,1)_{i,1}$ urn and $X_{0}=1$ to represent the initial
white ball in the urn. We will then define random variables
$M_{n}^{j,k}$ such
that
%
%
\begin{equation}
\law\bigl(M_{n}^{j,k} \bigr)=\law(W_n\mid
X_{j}=X_{k}=1), \label{50}
\end{equation}
so that by Lemma~\ref{lem20}, if $J$ and $K$ are random variables
independent of
$M_{n}^{j,k}$ satisfying
%
%
\begin{equation}
\label{51} \IP(K=k, J=j) = \frac{\IE(X_{k}X_{j})}{b^2}, \qquad\mbox
{$j,k\geq0$}
\end{equation}
for $b^2:=\IE W^2$, then $M^{J,K}_{n}$ has the double size bias distribution
of~$W$.

In order to generate a variable satisfying \eq{50} for $j<k$, we use the
following lemma that yields a method to construct an urn process having
the law
of the $(2,0;1,1)_{i,1}$ urn process up to time $n$ conditional on
$X_{k}=X_{j}=1$. This conditioned process follows the law of the
$(2,0;1,1,)_{i,3}$ urn process up to (and including) draw $j-1$. At
draw $j$,
exactly one \textit{black} ball is added and then draws $j+1$ through
$k-1$ follow the
$(2,0;1,1)$ urn law. Again at draw $k$ exactly one \textit{black} ball
is added
and then
the process continues to draw $n$ following the $(2,0;1,1)$ urn rule.
We write
$M_{n}^{j,k}$ to denote the number of white balls in the urn after $n$
draws in
this
process, and we refer to this process as the $M^{j,k}$ process. Our
next main
result shows
that this construction of $M_{n}^{j,k}$ has the distribution specified
in~\eq{50}. First we state a technical lemma; the proof can be found
at the end of this section.

%
\begin{lemma}\label{lem21}
Fix $i\geq1$ and let $W_n=\sum_{j=0}^n X_{j}$
where for $j\geq1$, $X_{j}$ is the indicator that a white ball is
drawn on
draw $j$ from the $(2,0;1,1)_{i,1}$ urn and $X_{0}=1$.
If $1\leq j <k \leq n$, then
%
%
\begin{equation}
\IP(X_j=1\mid X_k=1, W_{j-1}) =
\frac{1+W_{j-1}}{2j+i}, \label{52}
\end{equation}
and if $1\leq l < j <k \leq n$,
%
%
\begin{equation}
\IP(X_l=1\mid X_k=1, X_j=1,
W_{l-1}) = \frac{2+W_{l-1}}{2l+i+1}. \label{53}
\end{equation}
For all $1\leq j<k\leq n$,
\begin{eqnarray*}
\IE W_{n} &\leq&\sqrt{2\pi} \sqrt{\frac{n}{i+2}+
\frac{1}{2}},\qquad\IE X_{j} \leq\frac{\sqrt{\pi}}{\sqrt{(i+1)(i+2j-1)}},
\\
(2-\sqrt{\pi}) \frac{i+2n+1}{i+1}&\leq&\IE W_{n}^2 \leq2
\frac
{i+2n+1}{i+1},
\\
\IE(X_{j}X_{k}) &\leq& \frac{\sqrt{2\pi}(1+\sqrt{\pi})}{(i+2)\sqrt
{(i+2j)(i+2k-1)}}.
\end{eqnarray*}
If $\law(R_t)=(2,0;1,1)_{i,3}^t$, then for some constant $C$
independent of $t$
and $i$,
%
%
\begin{equation}
\IE R_t\leq C\sqrt{t/i}. \label{54}
\end{equation}
\end{lemma}
%
%
\begin{lemma}
Let $1\leq j<k \leq n$ and $M_n^{j,k}$, $W_n$ and $(X_l)_{l\geq1}$ be
defined as
in Lemma $\ref{lem21}$ and the preceding two paragraphs. Then
%
%
\[
\law\bigl(M_{n}^{j,k} \bigr)=\law(W_n\mid
X_{j}=X_{k}=1).
\]
\end{lemma}
\begin{pf}
Let $M_l^{j,k}=3+\sum_{t=1}^l m_t^{j,k}$, where
for $t\not=j,k$, $m_t^{j,k}$ is the indicator that draw $t$
in the $M^{j,k}$ urn process is white and $m_j^{j,k}=m_k^{j,k}=0$.
From the definition of the process,
for $l<j$,
%
%
\begin{equation}
\IP\bigl(m_{l}^{j,k}=1| M_{l-1}^{j,k}
\bigr)=\frac{M_l^{j,k}}{2l+i+1}. \label{55}
\end{equation}
%
And for $j<l<k$, since $m_j^{j,k}=0$,
%
%
\begin{equation}
\IP\bigl(m_{l}^{j,k}=1| M_{l-1}^{j,k}
\bigr)=\frac{M_l^{j,k}}{2l+i}. \label{56}
\end{equation}
%
Note also that $M_0^{j,k}=W_0+2=3$,
$m_j^{j,k}=m_k^{j,k}=0$
and draw $l>k$ in the $M^{j,k}$ urn process
follows the $(2,0;1,1)$ urn law.
Now comparing \eq{55} to \eq{53} and \eq{56} to \eq{52},
we find the sequential conditional probabilities agree
and so the lemma follows.
\end{pf}

We are now ready to prove Lemma~\ref{lem16}, and we first give the following
remark about the argument.
The $(2,0;1,1)_{i,3}$ process and the
$M^{j,k}$ process defined above differ only in that, in the latter process,
after each of draws $j$ and $k$ a single black ball is added into the urn
regardless of what is drawn; in the former process, the two balls added
to the
urn in these draws depend on the color drawn. This difference turns out
to be
small enough to allow a close coupling as stated in the lemma.

\begin{pf*}{Proof of Lemma~\ref{lem16}}
For each $t$ we construct
$ (r_{t}, m^{j,k}_{t} )$ to, respectively, denote the indicator
for the event that a white ball is added to the urn after draw number
$t$ for
the
$(2,0;1,1)_{i,3}$ process and for the $M^{j,k}$ process, and we write
%
%
\[
R_{n-1}=3+\sum_{t=1}^{n-1}
r_t, \qquad M_{n}^{j,k} =3+ \sum
_{t=1}^n m^{j,k}_{t}
\]
to denote the number of white balls in the urn after draw $n-1$ and $n$,
respectively, for each process.
Let $U_{t}$ be independent uniform
$(0,1)$ random variables. We define
%
%
\begin{equation}
\label{57} r_{t}=\I\biggl[U_{t}<\frac{R_{t-1}}{i+2t+1}
\biggr]
\end{equation}
and for $t\neq k, t\neq j$ we define
%
%
\[
m^{j,k}_{t} =\I\biggl[U_{t} <\frac{M^{j,k}_{t-1}}{i+2t+1-\I[t>j]-\I[t>k]}
\biggr].
\]
We also set $m^{j,k}_{k}=m^{j,k}_{j}=0$ since a single black ball is
added after
draws $j$ and $k$. Writing the event $M^{j,k}_{n} \neq{R}_{n-1}$ as a
union of
the events that index $t$ is the least
index such that $r_t\neq m^{j,k}_t$, and also using that
$m^{j,k}_{k}=m^{j,k}_{j}=0$,
we find that for $0<j<k$,
\begin{eqnarray*}
&&\IP\bigl(M^{j,k}_{n} \neq{R}_{n-1} \bigr) \\
&&\qquad\leq
\IE({r}_{k}+r_{j})
 +\IP\biggl(U_{n} < \frac{{R}_{n-1}}{i+2n-1} \biggr)\\
 &&\qquad\quad{} + \sum
_{t=j}^{n-1} \IP\biggl(\frac{R_{t-1}}{i+2t+1} <
U_{t} < \frac{{R}_{t-1}}{i+2t-1} \biggr).
\end{eqnarray*}
From \eq{57},
$\IE r_t = \IE R_{t-1} /(i+2t+1)$, so that we find
%
%
\begin{eqnarray*}
&&\IP\bigl(M^{j,k}_{n} \neq{R}_{n-1} \bigr)\\
&&\qquad\leq
\frac{\IE R_{k-1}}{i+2k+1}+\frac{\IE R_{j-1}}{i+2j+1}+\frac{\IE
R_{n-1}}{i+2n-1}
\\
&&\qquad\quad{}+ \sum_{t=j}^{n-1} \IE R_{t-1}
\biggl( \frac{1}{i+2t-1}-\frac{1}{i+2t+1} \biggr)
\\
&&\qquad\leq C \sqrt{\frac{k}{i}} \biggl(\frac{1}{i+2k+1} \biggr)+C \sqrt{
\frac{j}{i}} \biggl(\frac{1}{i+2j+1} \biggr)+C \sqrt{
\frac{n}{i}} \biggl(\frac{1}{i+2n-1} \biggr)
\\
&&\qquad\quad{}+ C \sum_{t=j}^{n-1} \sqrt{
\frac{t}{i}} \biggl( \frac{1}{i+2t-1}-\frac{1}{i+2t+1} \biggr)
\\
&&\qquad\leq C j^{-{ {1}/{2}}},
\end{eqnarray*}
where we have used \eq{54}.
Defining $J$ and $K$ as in \eq{51} we now have
%
%
\begin{eqnarray}\label{58}
&&\IP\bigl(M^{J,K}_{n}\neq{R}_{n-1} \bigr)
\nonumber
\\[-8pt]
\\[-8pt]
\nonumber
&&\qquad\leq\IP(J=0)+\IP(K=0)+ \IP(J = K)+\frac{2C}{b^2} \sum
_{j<k}j^{-{ {1}/{2}}}\IE X_jX_k.
\end{eqnarray}
Since $X_i^2=X_i$, $X_0=1$ and using \eq{51}, we have
\[
\IP(J=K)=\IP(K=0)=\IP(J=0)=\IE W_n /b^2,
\]
and now using the bounds on $\IE X_i X_j$, $\IE W_n$ and $\IE W_n^2$ from
Lemma~\ref{lem21},
we find that \eq{58} is bounded above by
\begin{eqnarray*}
\frac{3\IE W_n}{b^2}+\frac{C}{b^2} \sum_{j<k}
\frac{1}{\sqrt{jk}} \frac{1}{\sqrt{j}} \leq\frac{C}{\sqrt{n}}+\frac{C}{n}
\sum_{j<k} j^{-3/2} \leq\frac
{C}{\sqrt{n}}+C
\sum_{j} j^{-3/2} \leq\frac{C}{\sqrt{n}}.
\end{eqnarray*}
\upqed\end{pf*}

\begin{pf*}{Proof of Lemma \ref{lem21}} Let $A_m = \{X_{m}=1\}$. By
the definition of conditional
probability,
%
%
\begin{equation}
\label{59} \IP(A_l | A_k,A_j,W_{l-1})
= \frac{\IP(A_l\mid W_{l-1}) \IP(A_k A_j\mid A_l, W_{l-1})}{\IP(A_k A_j
\mid W_{l-1})}
\end{equation}
and
%
%
\begin{equation}
\label{60} \IP(A_j | A_k,W_{j-1}) =
\frac{\IP(A_j\mid W_{j-1}) \IP(A_k \mid A_j,
W_{j-1})}{\IP(A_k \mid W_{j-1})},
\end{equation}
and we next will calculate the probabilities above. For $j\geq1$, we have
%
%
\begin{equation}
\IP(A_j\mid W_{j-1}) = \frac{W_{j-1}}{i+2j-1}, \label{61}
\end{equation}
which implies that for $k\geq j$,
%
%
\begin{equation}
\label{62} \IP(A_k\mid W_{j-1}) = \frac{\IE(W_{k-1} \mid W_{j-1})}{i+2k-1}
\end{equation}
and
%
%
\begin{equation}
\label{63} \IP(A_k\mid A_j, W_{j-1}) =
\frac{\IE(W_{k-1} \mid A_j, W_{j-1})}{i+2k-1}.
\end{equation}
Now to compute the conditional expectations appearing above note first that
%
%
\begin{equation}
\label{64} \IE(W_{k}\mid W_{k-1}) =W_{k-1} +
\frac{W_{k-1}}{i+2k-1} = \biggl( \frac{i+2k}{i+2k-1} \biggr) W_{k-1}.
\end{equation}
Conditioning on $W_{j-1}$ and taking expectations yields
\[
\IE(W_{k}\mid W_{j-1}) = \biggl( \frac{i+2k}{i+2k-1}
\biggr)\IE(W_{k-1}| W_{j-1}),
\]
and then iterating and substituting $k-1$ for $k$ yields
%
%
\begin{equation}
\label{65} \IE(W_{k-1}\mid W_{j-1}) =\prod
_{t=1}^{k-j} \biggl( \frac{i+2(k-t)}{i+2(k-t)-1}
\biggr)W_{j-1}.
\end{equation}
Using \eq{65} with $j$ substituted for $j-1$, we also find
%
%
\begin{equation}
\label{66} \IE(W_{k-1}\mid A_j, W_{j-1}) =\prod
_{t=1}^{k-j-1} \biggl( \frac{i+2(k-t)}{i+2(k-t)-1}
\biggr) (1+W_{j-1});
\end{equation}
note here that conditioning on $A_j$ and $W_{j-1}$ is equivalent
to conditioning on $W_j$ and the event $\{W_j=W_{j-1}+1\}$.
We use a similar approach to obtain
\begin{eqnarray*}
\IE\bigl(W^2_{k}\mid W_{k-1}\bigr)
&=&W^2_{k-1} \biggl(1-\frac{W_{k-1}}{i+2k-1} \biggr) +
(W_{k-1}+1)^2\frac{W_{k-1}}{i+2k-1}
\\
&=& \biggl( \frac{i+2k+1}{i+2k-1} \biggr) W^2_{k-1}+
\frac{W_{k-1}}{i+2k-1},
\end{eqnarray*}
which can then be added to \eq{64} while letting $D_k = W_{k}(1+W_{k})$
to obtain
\[
\IE(D_k|W_{k-1}) = \frac{i+2k+1}{i+2k-1} D_{k-1},
\]
and thus
\[
\IE(D_k|W_{j-1}) = \frac{i+2k+1}{i+2k-1}
\IE(D_{k-1}| W_{j-1}).
\]
Iterating and substituting $k-1$ for $k$ gives
%
%
\begin{equation}\qquad
\label{67} \IE(D_{k-1}\mid W_{j-1}) =\frac{i+2k-1}{i+2j-1}
D_{j-1} = \frac{i+2k-1}{i+2j-1}W_{j-1}(1+W_{j-1}),
\end{equation}
and using \eq{67} with $j$ substituted for $j-1$, we also find
\[
\IE(D_{k-1}\mid A_j, W_{j-1}) =
\frac{i+2k-1}{i+2j+1} (W_{j-1}+1) (W_{j-1}+2).
\]
Letting
\[
c=\frac{1}{(i+2j-1)(i+2k-1)}\prod_{t=1}^{k-j-1}
\biggl( \frac{i+2(k-t)}{i+2(k-t)-1} \biggr)
\]
and applying \eq{61}, \eq{63}, \eq{66} and \eq{67} we have
%
%
\begin{eqnarray} \label{68}
\IP(A_jA_k \mid W_{l-1}) & =& \IE\bigl(
\IP(A_j|W_{j-1}) \IP(A_k|A_j,W_{j-1})|W_{l-1}
\bigr)
\nonumber
\\
& =& c\IE(D_{j-1}|W_{l-1})
\\
& =& c\frac{i+2j-1}{i+2l-1}W_{l-1}(1+W_{l-1}),\nonumber
\end{eqnarray}
and by substituting $l$ for $l-1$ in \eq{68}, we also find
%
%
\begin{equation}
\label{69} \IP(A_jA_k\mid A_l,
W_{l-1}) = c\frac{i+2j-1}{i+2l+1}(1+W_{l-1})
(2+W_{l-1}).
\end{equation}
Substituting \eq{61}--\eq{63}, \eq{65}, \eq{66}, \eq{68} and
\eq{69} appropriately into \eq{59} and \eq{60} proves \eq{52}
and \eq{53}.

From \eq{65} we have
%
%
\begin{equation}
\label{70} \IE W_{n}=\prod_{t=1}^n
\frac{i+2t}{i+2t-1} =\frac{\Gamma({(i+1)}/{2} )\Gamma(n+{
{(i+1)}/{2}}
+{ {1}/{2}} )}{\Gamma({ {(i+1)}/{2}}+{ {1}/{2}} )
\Gamma(n+{(i+1)}/{2} )},
\end{equation}
and from \eq{61} we find
\[
\IE X_{j}= \frac{\IE W_{j-1}}{i+2j-1}.
\]
Now using \eq{63} and \eq{66} yields
\[
\IE(X_{j}X_{k}) =\frac{\IE(1+W_{j-1})\IE X_{j}}{i+2k-1} \prod
_{t=j+1}^{k-1}\frac{i+2t}{i+2t-1},
\]
and using \eq{67} we find
%
%
\begin{equation}
\IE W_{n}^2 = 2 \frac{i+2n+1}{i+1}- \IE
W_{n}. \label{71}
\end{equation}
Lemma~\ref{lem4} applied to \eq{70} implies
%
%
\begin{equation}
\label{72} \frac{1}{\sqrt{\pi}}\sqrt{\frac{2n}{i+2}+1} \leq\prod
_{t=1}^n \frac{i+2t}{i+2t-1} \leq\sqrt{\pi}\sqrt{
\frac{2n}{i+2}+1},
\end{equation}
and collecting the appropriate facts above yields the bounds on $\IE W_n$,
$\IE W_n^2$, $\IE X_i X_j$ and $\IE X_i$.

For the bound on $\IE R_t$, an argument similar to \eq{65} leading to
\eq{70}
yields that
\begin{eqnarray*}
\IE R_t = 3 \prod_{m=1}^{t}
\frac{i+2m+2}{i+2m+1} =\frac{\Gamma({(i+3)}/{2} )\Gamma(t+{
{(i+3)}/{2}}
+{ {1}/{2}} )}{\Gamma({ {(i+3)}/{2}}+{ {1}/{2}} )
\Gamma(t+{(i+3)}/{2} )},
\end{eqnarray*}
which can be bounded by Lemma~\ref{lem4} resulting
in inequalities which are similar to~\eq{72}.
\end{pf*}

\section*{Acknowledgments}
We are grateful to Carl Morris for valuable advice on this project and an
anonymous referee for very helpful and detailed comments.
The authors have also learned that Robert Gaunt (University of Oxford) has
independently been working on a second order Stein operator for the
variance-gamma distribution, but there appears to be little overlap
between his work and ours.


%


\printaddresses

\end{document}